\documentclass{amsart}
\usepackage{amssymb,xspace}
\usepackage[figureposition=bottom,tableposition=bottom,skip=5pt]{caption}
\usepackage{graphicx,color} 
\usepackage[small]{diagrams}
\def\mat#1{\ensuremath{#1}\xspace}
\def\makedef#1#2{\expandafter\gdef\csname #1\endcsname{#2}}
\def\makelet#1#2{\expandafter\let\csname #1\expandafter\endcsname\csname #2\endcsname}
\def\defmath#1{\makelet{temp@#1}{#1}\makedef{#1}{\mat{\csname temp@#1\endcsname}}}

\def\defbb#1{\makedef{c#1}{\mat{\mathbb{#1}}}}		% usage \cA
\def\defcal#1{\makedef{l#1}{\mat{\mathcal{#1}}}}	% usage \lA
\def\bbcal#1{\defbb{#1}\defcal{#1}}
\bbcal A
\bbcal B
\bbcal C
\bbcal D
\bbcal E
\bbcal F
\bbcal G
\bbcal H
\bbcal I
\bbcal J
\bbcal K
\bbcal L
\bbcal M
\bbcal N
\bbcal O
\bbcal P
\bbcal Q
\bbcal R
\bbcal S
\bbcal T
\bbcal U
\bbcal V
\bbcal W
\bbcal X
\bbcal Y
\bbcal Z

\def\al{\mat{\alpha}}
\def\be{\mat{\beta}}
\def\ga{\mat{\gamma}}
\def\Ga{\mat{\Gamma}}

\def\De{\mat{\Delta}}

\def\hi{\mat{\chi}}

\def\La{\mat{\Lambda}}
\def\la{\mat{\lambda}}

\def\om{\mat{\omega}}

\def\si{\mat{\sigma}}
\def\Si{\mat{\Sigma}}

\def\te{\mat{\theta}}

\defmath{eta}
\defmath{xi}
\defmath{mu}
\defmath{nu}
\defmath{Phi}
\defmath{Psi}
\defmath{pi}
\defmath{rho}

\def\mrm@#1{\mat{\mathrm{#1}}}

\def\deffrak#1{\makedef{g#1}{\mat{\mathfrak{#1}}}} % usage \ga,\gA
\deffrak b
\deffrak c
\deffrak d
\deffrak D
\deffrak h
\deffrak m
\deffrak n
\deffrak p
\deffrak P
\deffrak q
\deffrak r
\deffrak s
\deffrak t
\deffrak u
\deffrak v

\def\DMO{\DeclareMathOperator}

\DMO{\Hom}{Hom}
\DMO{\RHom}{RHom}
\DMO{\lHom}{\lH\mathit{om}}
\DMO{\Ext}{Ext}
\DMO{\lExt}{\lE\mathit{xt}}
\DMO{\End}{End}
\DMO{\Aut}{Aut}
\DMO{\Fun}{Fun}
\DMO{\Tor}{Tor}
\DMO{\ext}{ext}

\DMO{\Ob}{Ob}
\DMO{\Mor}{Mor}
\DMO{\im}{im}
\DMO{\coim}{coim}
\DMO{\coker}{coker}
\DMO{\Arr}{Arr}
\DMO{\Id}{Id}
\DMO{\id}{id}
\DMO{\add}{add} % splitting of idempotents (karoubinization)
\DMO{\ind}{ind} % category of ind-objects
\DMO{\pro}{pro} % category of pro-objects
\DMO{\Map}{Map} %
\DMO{\Iso}{Iso} %
\DMO{\Isom}{Isom}%
\DMO{\Ind}{Ind}

\DMO{\Presh}{Presh}

\DMO\coalg{Coalg}
\DMO{\Rep}{Rep}
\DMO{\Cor}{Cor}

\DMO{\Mod}{Mod}
\DMO{\rad}{rad}
\DMO{\soc}{soc}
\DMO{\ann}{ann}

\DMO{\Spec}{Spec}
\DMO{\Specm}{Specm}
\DMO{\spec}{Spec}
\DMO{\Proj}{Proj}
\DMO{\supp}{supp}
\DMO{\Coh}{Coh}
\DMO{\coh}{coh}
\DMO{\Qcoh}{QCoh}
\DMO{\QCoh}{QCoh}
\DMO{\Pic}{Pic}
\DMO{\Div}{Div}
\DMO{\ch}{ch}
\DMO{\Hilb}{Hilb}
\DMO{\Fitt}{Fitt}
\DMO{\Quot}{Quot}

\DMO{\Gras}{Gr}
\DMO{\Flag}{Flag}

\DMO{\cone}{cone}
\DMO{\Tw}{Tw}
\DMO{\rank}{rk}
\DMO{\rk}{rk}
\DMO{\codim}{codim}
\DMO{\cov}{cov}
\DMO{\sgn}{sgn}
\DMO{\td}{td}
\DMO{\GL}{GL}
\DMO{\SL}{SL}

\DMO\Der{Der}
\DMO\der{Der}
\DMO\coder{Coder}
\DMO{\diag}{diag}
\DMO{\HMod}{HMod} %the homotopy category of modules over DGC
\DMO{\ad}{ad}
\DMO{\Ad}{Ad}
\DMO*{\colim}{colim}
\DMO*{\hocolim}{hocolim}
\DMO*{\holim}{holim}
\DMO{\Ho}{Ho}
\DMO{\har}{char}
\DMO{\sk}{sk}
\DMO{\cosk}{cosk}
\DMO{\Gal}{Gal}
\DMO{\tr}{tr}
\DMO{\Tr}{Tr}
\DMO{\Sh}{Sh}
\DMO{\Is}{Is} %Isometries
\DMO{\Hol}{Hol} %Holomorphic automorphisms
\DMO{\Lie}{Lie} %Lie algebra of a group
\DMO{\Res}{Res} %restriction, residue
\DMO{\irr}{irr} %
\DMO{\Irr}{Irr} %
\DMO{\Exp}{Exp} %
\DMO{\Log}{Log} %
\DMO{\Pow}{Pow}
\DMO{\pow}{pow}
\DMO{\mult}{mult} %
\DMO{\height}{ht} %
\DMO{\wt}{wt}
\DMO{\Vect}{Vect}
\DMO{\hd}{hd} %homological dimension
\DMO{\face}{face}

\newcommand{\GIT}{/\!\!/}

\def\dd{\mat{\partial}}
\def\iso{\simeq}
\def\ts{\otimes}

\def\sb{\subset}

\def\sp{\supset}

\def\xx{\times}

\def\ms{\backslash} %minus set
\DMO{\Eig}{Eig} 

\def\op{^{\rm op}}
\def\inv{^{-1}}
\def\dual{^\vee}

\def\ang#1{\mat{\left\langle #1\right\rangle}}

\def\set#1{\mat{\{ #1\}}}
\def\sets#1#2{\mat{\{ #1 \mid #2\}}}

\def\smat#1{\mat{\left(\begin{smallmatrix}#1\end{smallmatrix}\right)}}

\def\emb{\hookrightarrow}
\def\mto{\mapsto}
\def\arr{\futurelet\test\arrtest}
\def\arrtest{\ifx^\test\let\next\arra\else\let\next\arrb\fi\next}
\def\arra^#1{\xrightarrow{#1}} \def\arrb{\to}

\def\arrowsD{
\def\mto{{\:\vrule height .9ex depth -.2ex width .04em\!\!\!\;\ar}}
\def\ar{{\:\vrule depth -.52ex height .60ex width 0.85em\;\!\!\rhla\,}}
\def\arr{\futurelet\test\arrtest}
\def\arrtest{\ifx^\test\let\next\arra\else\let\next\arrb\fi\next}
\def\arra^##1{\rTo^{##1}} \def\arrb{\ar}
\def\emb{\futurelet\test\embtest}
\def\embtest{\ifx^\test\let\next\emba\else\let\next\embb\fi\next}
\def\emba^##1{\rInto^{##1}} \def\embb{{\:\rthooka\!\!\!\ar}}
\newarrow{Eq}=====
\def\rrarr{\pile{\rTo\\ \rTo}}
\def\lrarr{\pile{\rTo\\ \lTo}}   
\newarrow{ShortTo}{}{}-->
}

\def\arrowsDStandard{
\newarrow{TeXto}----{->}
\newarrow{TeXinto}C---{->}
\newarrow{TeXonto}----{->>}
\newarrow{TeXdashto}{}{dash}{}{dash}{->}
\newarrow{Eq}=====
\def\ar{\rightarrow}
\def\emb{\futurelet\test\embtest}
\def\embtest{\ifx^\test\let\next\emba\else\let\next\embb\fi\next}
\def\emba^##1{\rTeXinto^{##1}} \def\embb{\hookrightarrow}
\def\rrarr{\pile{\rTo\\ \rTo}}
\def\lrarr{\pile{\rTo\\ \lTo}}   
}

\newif\ifukr\ukrfalse
\newif\ifrus\rusfalse
\newif\ifger\gerfalse

\def\theorems{
\newcounter{nthr} %auxiliary counter
\numberwithin{nthr}{section}
\let\theHnthr\thenthr
\newtheorem{thr}[nthr]{Theorem}
\newtheorem{prp}[nthr]{Proposition}
\newtheorem{lmm}[nthr]{Lemma}
\newtheorem{crl}[nthr]{Corollary}
\newtheorem{clm}[nthr]{Claim}
\newtheorem{conj}[nthr]{Conjecture}
\theoremstyle{definition}
\newtheorem{dfn}[nthr]{Definition}
\newtheorem{rmr}[nthr]{Remark}
\newtheorem{exm}[nthr]{Example}
\newtheorem{claim}[nthr]{Claim}}

\def\ub#1{\mat{\overline{#1}}}  %upper bar
\def\lqq{\lq\lq}
\def\rqq{\rq\rq\xspace}

\def\eprint#1#2{% called from bbl file, #1 - archive name, #2 - preprint id
\expandafter\ifx\csname eprint@#1\endcsname\relax#1:#2\else%
\expandafter\def\csname\csname eprint@#1@idname\endcsname\endcsname{#2}%
\csname eprint@#1\endcsname\fi}

\def\defArchive#1#2#3{ %#1 - archive name, #2 - id name, #3 - code using id name
\makedef{eprint@#1}{#3}\makedef{eprint@#1@idname}{#2}}

\defArchive{arxiv}{id}{\href{http://arXiv.org/abs/\id}{{\tt arXiv:\id}}}
\defArchive{numdam}{id}{\href{http://www.numdam.org/item?id=\id}{Numdam:\id}}

\theorems
\usepackage{hyperref}

\graphicspath{{images/}}
{\makeatletter\gdef\input@path{{images/}}}

\DMO{\inn}{inn}
\DMO{\temp}{mod}\let\mod\temp
\def\a{A} % qp algebra
\def\b{B} % \ker(Z^{Q_0}\to Z)

\begin{document}
\title[Crepant resolutions and brane tilings]
{Crepant resolutions and brane tilings II: Tilting bundles}%

\author{Martin Bender}%
\author{Sergey Mozgovoy}%

%\address{}%

\email{mbender@wmaz.math.uni-wuppertal.de}%
\email{mozgov@math.uni-wuppertal.de}%
%\thanks{}%

%\subjclass[2000]{16G20}%
%\keywords{}%

%\date{16.08.2009}%
%\dedicatory{}%
%\commby{}%
% ----------------------------------------------------------------
\begin{abstract}
Given a brane tiling, that is, a bipartite graph on a torus,
we can associate with it a singular $3$-Calabi-Yau
variety.
Using the brane tiling, we can also construct all crepant
resolutions of the above variety.
We give an explicit toric description of tilting bundles
on these crepant resolutions. This result proves
the conjecture of Hanany, Herzog and Vegh and a version
of the conjecture of Aspinwall.
\end{abstract}
  \maketitle
%\tableofcontents
\section{Introduction}
The goal of this paper is to prove the conjecture
of Hanany, Herzog and Vegh \cite{HHV} on the description of
tilting bundles on the crepant resolutions
of singular $3$-Calabi-Yau varieties arising from brane tilings.
All these crepant resolutions can be constructed
as moduli spaces of representations of some
quiver with relations \cite[Theorem 15.1]{Ishii2}. 
These moduli spaces are toric $3$-Calabi-Yau varieties.
An explicit construction of their toric diagrams
was given in \cite{Mozgov6}.

Given a brane tiling, we can associate with it
a quiver potential $(Q,W)$ and a quiver potential
algebra $\cC Q/(\dd W)$.
The singular Calabi-Yau variety mentioned above
is isomorphic to the spectrum of the center of $\cC Q/(\dd W)$.
It has a non-commutative crepant
resolution $\cC Q/(\dd W)$ \cite{Bock2,Mozgov6}.
Its crepant resolutions are given by the moduli spaces
$\lM_\te=\lM_\te(\cC Q/(\dd W),\al)$
of $\te$-semistable $\cC Q/(\dd W)$-representations
of dimension $\al=(1,\dots,1)\in\cZ^{Q_0}$, where
$\te\in\cZ^{Q_0}$ is \al-generic.
All \te-semistable points in $\lM_\te$ are \te-stable for such \te.
Therefore there exists a universal
(also called tautological) vector bundle $\lU$ over $\lM_\te$,
endowed with a structure of a left $\cC Q/(\dd W)$-module.
It follows from the results of Van den Bergh (see \cite{Bergh1})
that $\lU$ is a tilting bundle (an alternative proof 
can be found in \cite{Ishii2}).
This vector
bundle can be decomposed into a sum of $\# Q_0$
line bundles. We will describe the toric Cartier
divisors inducing these line bundles. Namely,
we fix some $i_0\in Q_0$ and for every vertex $i\in Q_0$
we choose some path $u_i:i_0\arr i$. Intersecting
the path $u_i$ with perfect matchings
(they parametrize $2$-dimensional orbits of $\lM_\te$,
i.e.\ rays of the corresponding fan, see Section \ref{prelim}), we get a toric
Cartier divisor which induces a line bundle $\ub L_i$ over
$\lM_\te$. We will show that $\lU$ is isomorphic
to the direct sum of $\ub L_i$, $i\in Q_0$.
This description of the tilting bundle was conjectured
by Hanany, Herzog and Vegh \cite[Section~5.2]{HHV}.
Our result proves also a conjecture of Aspinwall \cite{Aspinwall}
on the existence of some \lqq globally defined\rqq
collection of line bundles that gives rise
to the tilting collection on $\lM_\te$ for arbitrary generic \te.
We should note that a similar description
of the exceptional collections in the context
of toric quiver varieties (these are moduli spaces 
of quiver representation for quiver without relations)
was given by Altmann and Hille \cite{Altmann-Hille}.

The paper is organized as follows: 
In Section \ref{prelim} we gather preliminary material
on brane tilings and the induced quiver potential algebras.
In Section  \ref{sec:toric} we recall some results
of Thaddeus \cite{Thaddeus1} about toric
quotients of toric varieties and prove some
facts on the descent of line bundles with respect
to such quotients.
In Section \ref{sec:tilting} we give a toric description
of the tilting bundle on $\lM_\te$.
In Section \ref{sec:examples} we give some explicit examples.

We would like to thank Markus Reineke for many useful
discussions.

\section{Preliminaries}\label{prelim}
Most of the content of this section can be found in \cite{Mozgov6}.
We briefly recall some material for the convenience of the  reader.

A brane tiling is a bipartite graph $G=(G_0^\pm,G_1)$
together with an embedding of the corresponding CW-complex 
into the real two-dimensional torus $T$ so that the complement
$T\ms G$ consists of simply-connected components.
The set of connected components of $T\ms G$ is denoted by 
$G_2$ and is called the set of faces of $G$.

With any brane tiling we can associate a quiver $Q=(Q_0,Q_1)$
embedded in a torus $T$ and a potential $W$ (linear combination of cycles in $Q$),
see \cite{MR2}. The set $Q_2$ of connected components
of $T\ms Q$ is called the set of faces of $Q$. 
The summands of $W$ are the cycles along the faces of $Q$
taken with appropriate signs.
With this data, we associate a quiver potential algebra
$\a=\cC Q/(\dd W)$, see \cite{MR2}.

For any arrow $a\in Q_1$, we define $s(a),t(a)\in Q_0$ to be its
source and target (also called tail and head).
Consider a complex of abelian groups
$$\cZ^{Q_2}\arr^{d_2}\cZ^{Q_1}\arr^{d_1}\cZ^{Q_0},$$
where $d_2(F)=\sum_{a\in F}a$, $F\in Q_2$ and $d_1(a)=t(a)-s(a)$
for any arrow $a\in Q$. 
%The sign in t(a)-s(a) is very important, as it says how the torus (C^*)^{Q_0}
%acts on the representation space C^{Q_1}.
Its homology groups are isomorphic to the homology
groups of the $2$-dimensional torus containing $Q$.
We define an abelian group \La by a cocartesian left upper square
of the following diagram
\begin{diagram}
\cZ^{Q_2}&\rTo^{d_2}&\cZ^{Q_1}&\rTo^{d_1}&\cZ^{Q_0}\\
\dTo&&\dTo^{\wt}&\ruDotsto_d\\
\cZ&\rTo^{\om_\La}&\La\\
\dDotsto^{\om_M}&\ruTo_i\\
M
\end{diagram}
where the left arrow of the square is given by $F\mto1,\ F\in Q_2$.
There exists a unique map $d:\La\arr\cZ^{Q_0}$ making the right
triangle commutative. Let $M=\ker(d)$. There exists a
unique map $\om_M:\cZ\arr M$ making the lower triangle commutative.
If $G$ has at least one perfect matching
then $\La$ is a free abelian group and the map
$\om_\La:\cZ\arr\La$ is injective (see \cite[Lemma 3.3]{MR2}).

We define a weak path in $Q$ to be a path consisting of arrows
of $Q$ and their inverses (for any arrow $a$, 
we identify $aa\inv$ and $a\inv a$ with trivial paths). 
For any weak path $u$, we define its content $|u|\in\cZ^{Q_1}$
by counting the arrows of $u$ with appropriate signs.
We define the weight of $u$
to be $\wt(u)=\wt(|u|)\in\La$. We define $\ub\om\in\La$ to
be the weight of any cycle along some face of $Q$. Note
that $\ub\om=\om_\La(1)$ and that $\ub\om\in M$.  
 
Let $\b=\ker(\cZ^{Q_0}\arr\cZ)$, where the map is given by
$i\mto1$, $i\in Q_0$. This group is generated by the elements
of the form $i-j$, where $i,j\in Q_0$.
As $Q$ is connected, we conclude that $\b=\im d$.
There is a short exact sequence
$$0\arr M\arr^i\La\arr^d \b\arr 0.$$
One can easily see that $\rk B=\#Q_0-1$, $\rk\La=\#Q_0+2$, and $\rk M=3$.

Let $\La^+\sb\La$ be a semigroup generated by the weights
of the arrows. Let $P\sb\La_\cQ$ be a cone generated by
$\La^+$
$$P=\sets{\sum a_i\la_i}{a_i\in \cQ_{\ge0},\ \la_i\in\La^+\text{ for all }i}.$$
We define $M^+=\La^+\cap M$ and $P_M=P\cap M_\cQ$.
We do not have $\La^+=P\cap\La$ in general.
But we have $M^+=P_M\cap M$ \cite{Mozgov6}.
This implies that $\Spec\cC[M^+]$ is a normal toric variety.
If the brane tiling is consistent (see e.g.\ \cite{Mozgov6})
then the quiver potential algebra $\cC Q/(\dd W)$ is
a $3$-Calabi-Yau algebra (see \cite{MR2,Broomhead1,Davison1})
and is a non-commutative crepant resolution of $\Spec\cC[M^+]$
(see \cite{Bock2,Mozgov6}). In this paper we will
always assume that our brane tiling is consistent.

Let \lA be the set of all perfect matchings of the bipartite graph $G$.
Any perfect matching $I\in\lA$ can be considered as
a subset of $Q_1$. We define a characteristic function
$\hi_I:\cZ^{Q_1}\arr\cZ$ by the rule (for $a\in Q_1$)
$$\hi_I(a)= 
\begin{cases}
	1,&a\in I,\\
	0,&a\not\in I.
\end{cases}$$
For any face $F\in Q_2$ we have $\hi_I(d_2(F))=1$.
Therefore we can consider $\hi_I$ as a linear map
$\La\arr\cZ$, i.e.\ as an element $\hi_I\in\La\dual$.
We define $\ub\hi_I=i^*\hi_I\in M\dual$.
The family of all $\hi_I:\La\arr\cZ$ (resp.\ 
$\ub\hi_i:M\arr\cZ$) defines a linear map $\hi:\La\arr\cZ^\lA$
(resp.\ $\ub\hi:M\arr\cZ^\lA$).

The dual cone $P\dual\in\La_\cQ\dual$
is generated by $\hi_I$, $I\in\La$,
and all the corresponding rays are extremal in $P\dual$
\cite[Lemma 2.3.4]{Broomhead1}.
Analogously, the dual cone $P_M\dual\in M_\cQ\dual$ is generated by
$\ub\hi_I$, $I\in\lA$. We denote by $\lA^e\sb\lA$ the set
of perfect matchings $I$, such that $\ub\hi_I$
generates an extremal ray in $P_M\dual$. These perfect
matchings are called extremal.

All the crepant resolutions of $\Spec\cC[M^+]$
can be described as moduli spaces of stable representations
of $A=\cC Q/(\dd W)$ in the sense of King \cite{King1}.
Let $\al=(1,\dots,1)\in\cZ^{Q_0}$ and let
$\te\in\b$ (i.e.\ $\te\in\cZ^{Q_0}$ is such that $\te\cdot\al=0$).
Any $A$-module $X$ can be described by the set of vector
spaces $(X_i)_{i\in Q_0}$ and linear
maps $X_a:X_i\arr X_j$ for arrows $a:i\arr j$.
We define $\dim X=(\dim X_i)_{i\in Q_0}\in\cZ^{Q_0}$.
An $A$-module $X$ of dimension \al is called
\te-semistable (resp.\ \te-stable) if for any proper $A$-submodule
$0\ne Y\sb X$ we have $\te\cdot \dim Y\ge0$ (resp.\ $\te\cdot\dim Y>0$).
We say that \te is $\al$-generic if for any $0<\be<\al$ we have
$\te\cdot\be\ne0$. In this case all \te-semistable modules
of dimension \al are automatically stable.
One can construct the moduli space $\lM_\te=\lM_\te(A,\al)$
of \te-semistable $A$-modules of dimension \al \cite{King1}.
It is shown in \cite{Mozgov6} that, for $\al$-generic $\te\in\b$,
this moduli space is a toric variety
(with a dense subtorus $T_M=\Hom_\cZ(M,\cC^*)$, see Section \ref{sec:toric})
$$\lM_\te=\Spec\cC[\La^+]\GIT_\te T_B=\Spec\cC[P\cap\La]\GIT_\te T_B,$$
where $T_B=\Hom_\cZ(B,\cC^*)$.
In this case $\lM_\te$ is smooth and is
a crepant resolution of $\Spec\cC[M^+]$ (see \cite{Ishii1,Mozgov6}).

An explicit description of a fan of $\lM_\te$ was given
in \cite{Mozgov6}. For any $A$-module
$X=((X_i)_{i\in Q_0},(X_a)_{a\in Q_1})$,
we define its cosupport $I_X=\sets{a\in Q_1}{X_a=0}$. 
It was shown in \cite{Mozgov6} that every $T_M$-orbit of $\lM_\te$
is uniquely determined by the cosupport of its modules.
Any such cosupport $I$ can be considered as
a subgraph of the bipartite graph $G$. 
It was shown in \cite{Mozgov6} that this subgraph can have
at most one connected component containing more than one edge
(we call it a big component of $I$).

\begin{prp}[{see \cite{Mozgov6}}]
Let $X\in\lM_\te$ and let $O_X\sb\lM_\te$ be its $T_M$-orbit.
Then
\begin{enumerate}
	\item $\dim O_X=3$ if and only if $I_X=\emptyset$.
	\item $\dim O_X=2$ if and only if $I_X$ is a perfect matching.
	\item $\dim O_X=1$ if and only if $I_X$ contains a big component which is a cycle.
	In this case $I_X$ is a union of two perfect matchings.
	\item $\dim O_X=0$ if and only if $I_X$ contains a big component which has two
	trivalent vertices of different colors and all other vertices
	of valence $2$. In this case $I_X$ is a union of three perfect
	matchings.
\end{enumerate}
\end{prp}

For any subset $I\sb Q_1$, we define a $\cC Q$-representation
$X_I=(X_{I,a})_{a\in Q_1}$
of dimension \al by the rule (for $a\in Q_1$)
$$X_{I,a}=
\begin{cases}
	0,&a\in I,\\
	1,&a\not\in I.
\end{cases}$$
We say that $I$ is $W$-compatible if $X_I$ is an $A$-representation.
For example, all perfect matchings and an empty set are $W$-compatible.
We say that $I$ is \te-stable
if $I_X$ is \te-stable.
The elements of the fan of $\lM_\te$ are in bijection
with $W$-compatible \te-stable subsets of $Q_1$.
The rays of the fan of $\lM_\te$
are parametrized by \te-stable perfect matchings.
All elements $\ub\hi_I\in M\dual$, $I\in\lA$, are contained in the
hyperplane
$$\sets{y\in M_\cQ\dual}{\om_M^*(y)=1},$$
where $\om_M:\cZ\arr M$ was defined earlier.
This implies that $\lM_\te$ is a toric
$3$-Calabi-Yau variety. The above proposition
gives an algorithm to construct its toric diagram
(this is an intersection of cones of the fan 
of $\lM_\te$ with the above hyperplane).
 %prelim
\section{Toric quotients} \label{sec:toric}
In this section we will recall some facts
from \cite{Thaddeus1} about toric quotients of
toric varieties and give further information
on the line bundles on such quotients.
We refer to \cite{Fulton2} and \cite{Oda1} for the relevant
definitions and properties of toric varieties.

Consider a pair $(\La,P)$, where \La is a lattice
(i.e.\ a free abelian group of finite rank) and
$P\sb\La_\cQ$ is a polyhedral cone.
We associate with it a scheme
$$X_P=X(\La,P):=\Spec \cC[P\cap\La].$$

More generally, given a pair $(\La,P)$, where \La
is a lattice and $P\sb\La_\cQ$ is a polyhedron, we associate
with it a scheme $X(\La,P)$ in the following way.
Let $C(P)\sb\cQ\xx\La_\cQ$ be a cone which is a closure of
$$\sets{\la(1,x)}{\la\in\cQ_{\ge0},x\in P}.$$
We endow $\cC[C(P)\cap(\cZ\xx\La)]$ with a \cZ-grading induced
by the first coordinate and define
$$X_P=X(\La,P):=\Proj \cC[C(P)\cap(\cZ\xx\La)].$$
Let $T_\La=\Hom_\cZ(\La,\cC^*)$.
There is a canonical $T_\La$-action on $X_P$ and
a canonical $T_\La$-linearization of the canonical line bundle
$\lO(1)$ on $X_P$.
If $P$ is a cone then $C(P)=\cQ_{\ge0}\xx P$ and
$C(P)\cap(\cZ\xx\La)=\cN\xx(P\cap\La)$. So our new
definition of $X_P$ is compatible with the old
one.

The scheme $X_P$ can be described  as a toric variety 
associated to a fan.
For any $y\in \La\dual_\cQ$ define
$$\face_y(P):=\sets{x\in P}{\ang{x,y}=\min_P\ang{-,y}}.$$
All faces of $P$ have this form for some $y\in\La\dual_\cQ$.
For any face $F\sb P$, define its normal cone 
$$N_F=N_FP:=\sets{y\in\La\dual_\cQ}{\face_y(P)\sp F}=\sets{y\in\La\dual_\cQ}{\ang{F,y}\le\ang{P,y}}.$$
For any faces $F,G\sb P$, we have $F\sb G$ if and only if $N_GP\sb N_FP$.
The set of cones 
$$N(P)=\sets{N_FP}{F\text{ face of }P}$$
is a fan in $\La\dual$ and the associated toric
variety is isomorphic to $X_P$ 
(cf.\ \cite[Prop.\ 2.17]{Thaddeus1}).

\begin{lmm}
Let $F\sb P$ be some face and let
$\ang F\sb\La_\cQ$ be a vector space generated by 
the differences $x-y$, for $x,y\in F$.
Then $\ang F= N_F^\perp$.
\end{lmm}
\begin{proof}
We can suppose that $0\in F$. Then $y\in N_F$ if and only if
$$\ang{F,y}=0,\qquad \ang{P,y}\ge0.$$
This implies that $\ang F\sb N_F^\perp$.
This vetor spaces have equal dimension as $\dim F+\dim N_F=\dim\La_\cQ$.
\end{proof}

For any face $F\sb P$, let $O_F$ denote the $T_\La$-orbit
corresponding to $N_F$.

\begin{lmm}[see {\cite[Prop.\ 2.13]{Thaddeus1})}]
Let $F\sb P$ be some face. Then
\begin{enumerate}
	\item $\dim O_F=\dim F$.
	\item The character group of the stabilizer of $O_F$ in $T_\La$
	equals $\coker(N_F^\perp\cap\La\arr\La)$.
	\item The closure of $O_F$ equals $X(\La,F)$.
	\item For any two faces $F,G$ we have $F\sb G$ if and only if $O_F\sb O_G$.
\end{enumerate}
\end{lmm}

Consider an exact sequence of lattices
$$0\arr M\arr^i \La\arr^d \b\arr 0.$$
For any face $F\sb P$, we define $F_M=F\cap M_\cQ$.
There is an inclusion $T_\b\sb T_\La$
that induces an action of $T_\b$ on $X_P$ and
on the line bundle $\lO(1)$.
It is shown in \cite[Prop.\ 3.2]{Thaddeus1} that
the corresponding GIT quotient is given by
$$X(\La,P)\GIT T_\b=X(M,P_M).$$

\begin{lmm}[{\cite[Lemma 3.3]{Thaddeus1}}]
Let $F\sb P$ be a face. Then
\begin{enumerate}
	\item $O_F$ is $T_\b$-semistable if and only if $F\cap M_\cQ\ne\emptyset$.
	\item $O_F$ is $T_\b$-stable if and only if $M_\cQ$ intersects $\inn(F)$
	transversally.
	\item The image of a $T_\b$-semistable orbit $O_F$ in 
	$X_P\GIT T_\b$ is $O_{F_M}$.
\end{enumerate}
\end{lmm}

\begin{rmr}
Condition that $M_\cQ$ intersects $\inn(F)$ transversally
means that\\
$\inn(F)\cap M_\cQ\ne\emptyset$ and $\ang{F}+M_\cQ=\La_\cQ$.
We say that $F$ is $M$-stable in this case.
We denote the subscheme of stable points of $X_P$
by $X^s_P$.
\end{rmr}

\begin{rmr}
There is a bijection between the faces of $P_M$ and
the faces $F\sb P$ such that $\inn(F)\cap M_\cQ\ne\emptyset$.
\end{rmr}

\begin{prp}
The set of cones
$$N^s(P_M)=\sets{N_{F_M}P_M\sb M\dual_\cQ}{F\sb P\text{ is }M\text{-stable}}$$
forms a fan in $M\dual$. The corresponding toric variety
is isomorphic to $X^s_P\GIT T_\b$.
\end{prp}
\begin{proof}
The set $N^s(P_M)$ is a subset of the fan $N(P_M)$.
To show that $N^s(P_M)$ is a fan, we just have to show
that any face of the cone from $N^s(P_M)$ is contained
in $N^s(P_M)$.
Let $F\sb P$ be $M$-stable. Let $\tau'\sb N_{F_M}P_M$
be some face. We will show that $\tau'\in N^s(P_M)$.
We can find some face $G'\sb P_M$ such that $\tau'=N_{G'}P_M$.
We choose a minimal face $G\sb P$ such that $G_M=G'$.
The minimality property implies that $\inn{G}\cap M_\cQ\ne\emptyset$.
Moreover, $N_{G_M}P_M\sb N_{F_M}P_M$, so $F_M\sb G_M$ and therefore $F\sb G$.
In particular, $\ang G+M_\cQ=\La_\cQ$ and therefore $G$ is $M$-stable.
This implies that $\tau'\sb N^s(P_M)$.
\end{proof}

\begin{lmm}\label{lmm:bijection}
Let $F\sb P$ be an $M$-stable face.
Consider the normal cones $N_F=N_FP\sb\La\dual_\cQ$ 
and $N_{F_M}=N_{F_M}P_M\sb M\dual_\cQ$.
Then the map $i^*:\La\dual_\cQ\arr M\dual_\cQ$ restricts to
a bijection
$$i^*_F:N_F\arr N_{F_M}.$$
\end{lmm}
\begin{proof}
Without loss of generality we may assume that $0\in \inn(F)$.
Then for any $y\in N_F$, we have
$$\ang{F,y}=0,\qquad \ang{P,y}\ge0.$$
This implies that
$$\ang{F_M,i^*(y)}=0,\qquad \ang{P_M,i^*(y)}\ge0$$
and therefore $i^*(y)\in N_{F_M}$.

It follows from our assumption that the vector space 
$\ang F$ intersects $M_\cQ$ transversally.
This implies that the homomorphism of vector spaces
$$F^\perp\arr F_M^\perp$$
is an isomorphism. Therefore the 
map $i^*:N_F\arr N_{F_M}$ is injective,
as $N_F\sb F^\perp$.

Let us prove the surjectivity. Consider
$y'\in N_{F_M}\sb F_M^\perp$. We can find
$y\in F^\perp$ such that $i^*(y)=y'$.
We know that $\ang{P_M,y}\ge0$
and we have to show that $\ang{P,y}\ge0$, as this
will imply $y\in N_F$.
Assume that there exists $x_0\in P$ such that $\ang{x_0,y}<0$.
Without loss of generality we may assume that $\La_\cQ$
is generated by $x_0$ and $F$, and that $P$ is a convex hull
of $x_0$ and $F$. Then $y^\perp=\ang F$.
It follows from the transversality of the intersection
of $M_\cQ$ and $F$ that $M_\cQ$ intersects $P\ms F$.
But for any point $x$ in this intersection we have 
$\ang{x,y}<0$ and $x\in P_M$. This contradicts
our assumption $\ang{P_M,y}\ge0$.
\end{proof}

\begin{crl}
For any $M$-stable face $F\sb P$, we have
$$N_{F_M}\dual=N_F\dual\cap M_\cQ.$$
\end{crl}
\begin{proof}
We have
$$N_F\dual\cap M_\cQ=\sets{x\in M_\cQ}{\ang{x,N_F}\ge0}
=\sets{x\in M_\cQ}{\ang{x,N_{F_M}}\ge0}=N_{F_M}\dual.$$
\end{proof}

\begin{crl}
Let $F\sb P$ be an $M$-stable face. 
Let $\si=N_FP$, $\si'=N_{F_M}P_M$,
$U_\si=\cC[\si\dual\cap\La]$, and $U_{\si'}=\cC[(\si')\dual\cap M]$.
Then the map $X_P^{s}\arr X_P^s\GIT T_\b$
is given over $U_{\si'}$ by
$$U_{\si}=\Spec\cC[N_F\dual\cap \La]\arr \Spec\cC[N_F\dual\cap M]=U_{\si'}.$$
\end{crl}

Recall that with any $N(P)$-linear support function
$h:|N(P)|\arr\cQ$ (see e.g.\ \cite[Section 2.1]{Oda1})
we can associate a $T_\La$-equivariant
line bundle $L_h$ over $X_P$ (see \cite[Prop.\ 2.1]{Oda1}).
If $T_\b$ acts freely on $X^s_P$ then this line bundle descends
to a $T_M$-equivariant line bundle on
$X_P^s\GIT T_\b$
(this follows from \cite[Prop.~0.9]{GIT}
and descent theory). We are going to describe an
$N^s(P_M)$-linear support function that
gives this line bundle.

\begin{thr}\label{descend}
Assume that $T_\b$ acts freely on $X^s_P$.
Let $h:|N(P)|\arr\cQ$ be an $N(P)$-linear support function.
Define an $N^s(P_M)$-linear support function
$h':|N^s(P_M)|\arr\cQ$ by the rule
$$h'(y)=h((i^*_F)\inv(y)),\qquad y\in N_{F_M},$$
where $F\sb P$ is an $M$-stable face
and $i^*_F:N_F\arr N_{F_M}$ is a bijection defined earlier.
Then the descend of $L_h$ to $X^s_P\GIT T_\b$ is isomorphic to $L_{h'}$ as a
$T_M$-equivariant line bundle.
\end{thr}
\begin{proof}
Let us recall the construction of a $T_\La$-equivariant line
bundle $L_h$ on $X_P$ associated to the support function $h:|N(P)|\arr\cQ$
(see \cite[Prop.\ 2.1]{Oda1}).
For any commutative semigroup $\Ga$, we denote the canonical basis
of the semigroup algebra $\cC[\Ga]$ by $(e^\ga)_{\ga\in\Ga}$.
We can find a system of elements $(l_\si\in\La)_{\si\in N(P)}$
such that $h|_\si=l_\si|_\si$ for any $\si\in N(P)$.
The line bundle $L_h$ is defined by gluing the line bundles
$U_\si\xx\cC$ over $U_\si$, $\si\in N(P)$ using the gluing functions
$$g_{\tau\si}:(U_\si\xx\cC)|_{U_\tau}\arr U_\tau\xx\cC,\quad
(x,c)\mto(x,e^{l_\si-l_\tau}(x)c)$$
for $\tau<\si$.
The action of $T_\La$ on $U_\si\xx\cC$ is given by
$$t(x,c)=(tx,e^{-l_\si}(t)c),\quad t\in T_\La.$$
Let now $\si=N_F$ and $\si'=N_{F_M}$,
for some $M$-stable face $F\sb P$.
Let $\pi:U_\si\arr U_{\si'}=U_{\si}\GIT T_\b$ be a canonical projection
We give an explicit description of the descend
line bundle $(U_\si\xx\cC)\GIT T_\b$ over
$U_\si\GIT T_\b=U_{\si'}$.

The character group of the stabilizer of $O_F$ in $T_\b$
is given by $\coker(\si^\perp\cap\La\arr \b)$ 
(see e.g.\ \cite[Prop.\ 2.6]{Thaddeus1}). By our assumptions
this stabilizer is trivial, so 
$$(\si^\perp\cap \La)+M=\La.$$
This means that we can find some $m_\si\in M$ such that
$l_\si-m_\si\in\si^\perp$.
Consider the action of $T_M$ on 
$U_{\si'}\xx\cC$ given by
$$t(x,c)=(tx,e^{-m_\si}(t)c),\quad t\in T_M.$$
The map
$$\ub\pi:U_\si\xx\cC\arr U_{\si'}\xx\cC,\qquad (x,c)\mto (\pi(x),e^{l_\si-m_\si}(x)c),$$
is $T_\La$-equivariant. Indeed, for any $t\in T_\La$ we have
$$\ub\pi(t(x,c))=\ub\pi(tx,e^{-l_\si}(t)c)=(\pi(x),e^{l_\si-m_\si}(tx)e^{-l_\si}(t)c)
=(\pi(x),e^{-m_\si}(t)e^{l_\si-m_\si}(x)c).$$
On the other hand 
$$t\ub\pi((x,c))=t(\pi(x),e^{l_\si-m_\si}(x)c)=(\pi(x),e^{-m_\si}(t)e^{l_\si-m_\si}(x)c).$$
This shows that $\ub\pi:U_\si\xx\cC\arr U_{\si'}\xx\cC$ is
a quotient with respect to the action of $T_\b$.

The gluing of line bundles $(U_\si\GIT T_\b)\xx\cC$
is induced by the gluing of line bundles $U_\si\xx\cC$
and is given by the formula
$$U_{\tau'}\xx\cC\arr(U_{\si'}\xx\cC)|_{U_{\tau'}},\quad
(x,c)\mto(x,e^{m_\si-m_\tau}(x)c),$$
where $\si=N_F$, $\si'=N_{F_M}$, $\tau=N_G$, $\tau'=N_{G_M}$
for $M$-stable faces $F\sb G$ of $P$.
The corresponding support function $h':|N(P_M)|\arr\cQ$
is given on $y'\in\si'$ by
$$h'(y')=m_\si(y')=m_\si((i^*_F)\inv(y'))=l_\si((i^*_F)\inv(y'))=h((i^*_F)\inv(y')),$$
as $l_\si-m_\si\in\si^\perp$.
\end{proof}

%\begin{crl}
%Assume that $P$ is $M$-generic.
%Then there is a bijection between the faces of $P$ intersecting $M_\cQ$
%and the faces of $P_M$ given by
%$$F\mto F_M:=F\cap M_\cQ.$$
%\end{crl}
%\begin{proof}
%For any face $F'$ of $P_M$, the corresponding $T_M$-orbit $O_{F'}$ is the image
%of some stable orbit $O_F$ in $X_P$. This image is given by $O_{F\cap M_\cQ}$,
%so $F'=F\cap M_\cQ$. This proves the surjectivity of
%the above map. The injectivity follows from the fact that
%if $M_\cQ$ intersects face $F$ then it intersects $\inn(F)$. 
%\end{proof}

Any element $\te\in \b$ can be considered as a
character $\te:T_\b\arr\cC^*$. We can tensor the action
of $T_\b$ on $\lO(1)$ with this character.
The stable (resp.\ semistable) points of $X_P$ with
respect to this linearization are called \te-stable 
(resp.\ \te-semistable).
The corresponding GIT quotient is denoted by
$X_P\GIT_\te T_\b$. We have (see \cite[2.16]{Thaddeus1})
$$X_P\GIT_\te T_\b\iso X(\La,P^\te)\GIT T_\b=
X(M,P^\te\cap M_\cQ),$$
where $P^\te=P-\la$ for some $\la\in\La$ with $d(\la)=\te$.

%\begin{dfn}
%An element $\te\in \b$ is called $P$-generic if
%all \te-semistable points of $X_P$ are \te-stable.
%Polyhedron $P$ is called $M$-generic if all semistable
%(with respect to the canonical linearization of $\lO(1)$)
%points of $X_P$ are stable.
%\end{dfn}
 %toric quotients
\section{Tilting bundles}\label{sec:tilting}
Let $(Q,W)$ be a quiver potential
associated to some consistent brane tiling,
let $\a=\cC Q/(\dd W)$, and let $\al=(1,\dots,1)\in\cZ^{Q_0}$.
The goal of this section is to give a toric description
of the tilting bundles on the moduli spaces $\lM_\te(A,\al)$.

\begin{dfn}
Let $X$ be an algebraic variety.
A coherent sheaf $T\in\Coh X$ is called a tilting sheaf
if $\Ext^n(T,T)=0$ for $n>0$ and the triangulated
category $D^b(X)=D^b(\Coh X)$ is generated by
the summands of $T$. A collection of coherent sheaves 
$(T_i)_{i\in I}$ is called a tilting collection
if $\Ext^n(T_i,T_j)=0$ for $n>0$ and $i,j\in I$,
and the triangulated category $D^b(X)$ is
generated by the objects $T_i$, $i\in I$.
\end{dfn}

We use notation from Section \ref{prelim}.
In particular, we have defined
an exact sequence of free abelian groups
$$0\arr M\arr^i\La\arr^d\b\arr0$$
and a cone $P\sb\La_\cQ$ there.
Let $\te\in \b$ be \al-generic.
We have seen that the moduli space
$\lM_\te=\lM_\te(\a,\al)$ is a toric quotient
$$\lM_\te=\Spec\cC[P\cap\La]\GIT_\te T_\b
=X_P\GIT_{\te}T_\b=X(\La,P^\te)\GIT T_\b=X(M,P^\te\cap M_\cQ),$$
where $P^\te=P-\la$ for some
$\la\in\La$ with $d(\la)=\te$.

We know already how to parametrize the set of $T_M$-orbits of
$\lM_\te$, or equivalently, the fan of $\lM_\te$.
The set of rays of the fan of $\lM_\te$
is in bijection with \te-stable perfect matchings.
It is also in bijection
with the the set of facets (codimension $1$ faces)
of $P^\te_M=P^\te\cap M_\cQ$.

For any weak path $u$ and for any perfect matching $I\in\lA$,
we define $\hi_I(u)=\hi_I(\wt(u))$.
The extremal rays of $P\dual$ are parametrized by
the perfect matchings (see \cite[Lemma 2.3.4]{Broomhead1}).
This implies that for any weak path $u$ the system of integers
$(\hi_I(u))_{I\in\lA}$ determines a $T_\La$-Cartier divisor 
and therefore a $T_\La$-equivariant line bundle over $X_P$
which we denote by $L(u)$
(forgetting the $T_\La$-action, we get just a trivial line bundle).
If we restrict this system of integers to \te-stable perfect matchings,
we get a $T_M$-Cartier divisor and a $T_M$-equivariant
line bundle over $\lM_\te$ which we denote by $\ub L(u)$.

Let us fix some vertex $i_0\in Q_0$. For any
vertex $i\in Q$ we fix some weak path 
$u_i:i_0\arr i$.
The following result proves a conjecture of Hanany,
Herzog and Vegh \cite[Section~5.2]{HHV}

\begin{thr}
For any \al-generic $\te\in\b$,
the line bundles $\ub L(u_i)$, $i\in Q_0$, form
a tilting collection on $\lM_\te(A,\al)$.
\end{thr}
\begin{proof}
We know from \cite[Theorem 6.3.1]{Bergh1} 
and the fact that $A=\cC Q/(\dd W)$ is a non-commutative
crepant resolution of its center \cite{Mozgov6,Bock2}
that there is an equivalence of categories
$$\Psi:D^b(\mod \a\op)\arr D^b(\Coh \lM_\te),\quad M\mto M\ts_\a^L\lU,$$
where $\lU$ is a universal vector bundle on $\lM_\te$ (see also \cite{Ishii2}).
This implies, in particular, that the vector bundle $\lU=\Psi(\a)$
is a tilting sheaf.
We will give its toric description.
Let us recall the construction of the universal vector
bundle from \cite[Prop.~5.3]{King1}.

Let $(e_i)_{i\in Q_0}$ be the canonical basis of $\cZ^{Q_0}$
and let $T_0:=\Hom(\cZ^{Q_0},\cC^*)=\GL_\al(\cC)$.
Let $R=R(\a,\al)$ and let $R^s\sb R$ be the subvariety of
\te-stable representations. The diagonal $\De=\cC^*\sb T_0$
acts trivially on $R$.
We have $T_\b=T_0/\De$ and $\lM_\te=R\GIT_\te T_\b=R^s/ T_\b$.

For any $i\in Q_0$, we define a $T_0$-equivariant line bundle
$L_i$ over $R$ to be $R\xx\cC$ with an action on the second
factor induced by $e_i$. Explicitly, the action is given by
$$t(x,c)=(tx,t_ic),\qquad t=(t_i)_{i\in Q_0}\in T_0,\ (x,c)\in R\xx\cC.$$
The action of $T_\b$ on $L_i$ is not well defined as $\De$ acts
nontrivially on $L_i$.
Namely, it acts with weight $1$ on the second factor.
To overcome this problem,
we can multiply the action of $T_0$ with an
action going through 
some homomorphism $T_0\arr\De$ such that the
new action restricted to \De is trivial 
(see \cite[Prop.~5.3]{King1}, note that the $T_0$-orbits
will not change).
The homomorphism $\psi:T_0\arr\De$ is a character of $T_0$, that is,
an element $\psi\in\cZ^{Q_0}$. The triviality
condition means that $\psi\cdot\al=-1$. 
We make the choice $\psi=e_{i_0}$.
Then the action of $T_\b=T_0/\De$ on $L_i$ is given
by the character $e_i-e_{i_0}\in \b$.
Let the $T_\b$-equivariant line bundle $L_i$ on $R$
descend to the line bundle $\ub L_i$ on $\lM_\te=R^s/ T_\b$.
It is shown in \cite[Prop.~5.3]{King1} that 
$\oplus_{i\in Q_0}\ub L_i$ is a universal vector bundle on $\lM_\te$.

There is a natural action of $T_\La$ on $R$.
In order to extend it to an action on $L_i=R\xx\cC$
compatible with an action of $T_\b$,
we have to choose some $\la_i\in\La$ 
such that $d(\la)=e_i-e_0$.
We choose $\la_i=\wt(u_i)\in\La$.
The inverse image of $L_i$ with respect to the
natural map $X_P\arr R$ 
(this is a normalization of some irreducible component of $R$,
see \cite{Mozgov6}) is given by $L(u_i)$.
This implies that the descent line bundle of $L_i$ with respect
to $R^s\arr\lM_\te$ is isomorphic to the descent line bundle of
$L(u_i)$ with respect to $X_P^s\arr R^s\arr\lM_\te$.
According to Theorem \ref{descend}, the 
descent line bundle of $L(u_i)$ is $\ub L(u_i)$.
This means that $\ub L_i\iso\ub L(u_i)$.
Therefore $\lU\iso\oplus_{i\in Q_0}\ub L(u_i)$.
\end{proof}

\begin{rmr}
If $u,v:i\arr j$ are two weak paths then $uv\inv$ is a weak cycle.
This implies that $\wt(u)-\wt(v)\in M$ and therefore
$\ub L(u)$ and $\ub L(v)$ are isomorphic line bundles
(see \cite[Section~3.4]{Fulton2}).
If we substitute the vertex $i_0$ by some vertex $i'_0$,
then the line bundles $L(u_i),\ i\in Q_0$ should be
tensored with a line bundle $L(u)$, where $u:i'_0\arr i_0$
is any weak path. This ambiguity corresponds to the ambiguity
of the universal vector bundle over $\lM_\te$. The universal
vector bundle is defined only up to tensoring with a line
bundle.
\end{rmr}

\begin{rmr}
The conjecture of \cite[Section~5.2]{HHV}
states actually that the collection $\ub L(u_i)$, $i\in Q_0$,
is an exceptional collection.
But this is certainly false, as 
$\Hom_{\lM_\te}(\ub L(u_i),\ub L(u_j))=e_j \a e_i$ (see Corollary \ref{crl:ext})
is always nonzero.
\end{rmr}

\begin{rmr}
The collection of line bundles $L(u_i),\ i\in Q_0$ on $X_P$ has a property
that for any \al-generic $\te\in \b$ it descends to a tilting
collection on $\lM_\te$. The existence of such \lqq globally defined\rqq
collection was conjectured by Aspinwall \cite{Aspinwall}.
\end{rmr}

\begin{crl}\label{crl:ext}
For any weak path $u:i\arr j$ in $Q$, we have 
\begin{enumerate}
	\item $H^n(\lM_\te,\ub L(u))=0$, $n>0$.
	\item $H^0(\lM_\te,\ub L(u))=e_j\a e_i$, where $\a=\cQ/(\dd W)$.
\end{enumerate}
\end{crl}
\begin{proof}
Let $\a=\cC Q/(\dd W)$.
By the proof of the above theorem, the vector bundle 
$$\lU=\oplus_{k\in Q_0}\ub L(u_k)$$
can be endowed with a structure of a universal vector bundle.
The map $\Psi:D^b(\mod \a\op)\arr D^b(\Coh \lM_\te)$
maps the right $\a$-module $e_k\a$ to the summand $\ub L(u_k)$ of $\lU$.
This implies, for $n>0$, 
\begin{multline*}
\Ext^n_{\a\op}(e_i\a,e_j\a)=0=\Ext^n_{\lM_\te}(\ub L(u_i),\ub L(u_j))\\
=\Ext^n_{\lM_\te}(\lO,\ub L(u_ju_i\inv))=H^n(\lM_\te,\ub L(u)).
\end{multline*}
For the $\Hom$-space we get
\begin{multline*}
\Hom_{\a\op}(e_i\a,e_j\a)=e_j\a e_i=\Hom_{\lM_\te}(\ub L(u_i),\ub L(u_j))\\
=\Hom_{\lM_\te}(\lO,\ub L(u_ju_i\inv))=H^0(\lM_\te,\ub L(u)).
\end{multline*}
\end{proof}

 %exceptional sequences
\section{Examples}\label{sec:examples}
In this section we will consider two examples:
the suspended pinch point and the quotient
singularity $\cC^3/(\cZ_2\xx\cZ_2)$.
In the first example we will study all possible
generic stabilities and in the second
example we will study only three particular stabilities.

\subsection{Suspended pinch point}
\begin{figure}[h!]%
\scalebox{0.8}{\includegraphics{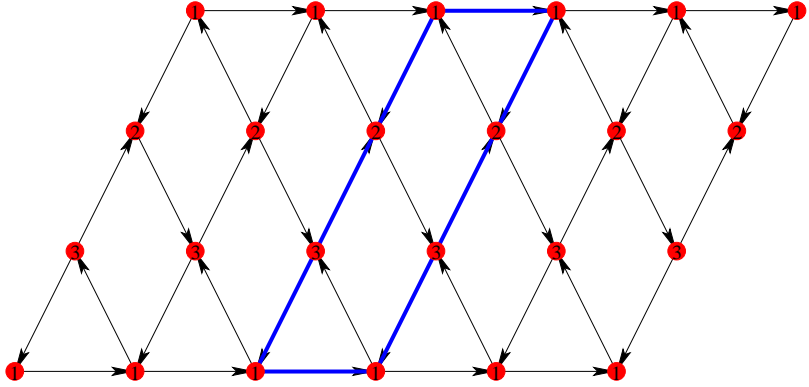}}%
\caption{Periodic quiver and a fundamental domain for SPP}%
\label{fig:spp}%
\end{figure}

Here we consider a brane tiling 
called a suspended pinch point. The
corresponding periodic quiver with a fundamental
domain is given in Figure \ref{fig:spp}.

%The underlying quiver is
%\begin{center}
%\scalebox{0.6}{\includegraphics{qspp}}
%\end{center}

Let $Q=(Q_0,Q_1,Q_2)$ be the corresponding quiver embedded
in a torus. We will denote the arrow from a vertex $i\in Q_0$
to a vertex $j\in Q_0$ by $ij$.
The list of all perfect matchings of the brane tiling
is given in Table \ref{pm:spp}. Every perfect matching
is described there as a subset of $Q_1$.

\begin{table}[h!]
\begin{tabular}{|c|c|}
\hline
$N$ &$I$\\
\hline
$1$&$12,31$\\
$2$&$21,13$\\
$3$&$32,11$\\
$4$&$23,11$\\
$5$&$12,13$\\
$6$&$21,31$\\
\hline
\end{tabular}
\caption{Perfect matchings of SPP}
\label{pm:spp}
\end{table}

Recall that we have defined a linear map $\ub\hi:M\arr\cZ^\lA$
in Section \ref{prelim}.
We can choose such basis of $M$ that the map $\ub\hi^*:\cZ^\lA\arr M\dual$
is given by the matrix
$$\smat{
0&2&1&0&1&1\\
0&0&1&1&0&0\\
1&1&1&1&1&1}$$
and the map $\om_M^*:M\dual\arr\cZ$ is given by the matrix 
$\smat{0&0&1}$. This gives us the toric diagram of $\Spec\cC[M^+]$
\begin{center}
\scalebox{1}{\begin{picture}(0,0)%
\includegraphics{spp.toric.pstex}%
\end{picture}%
\setlength{\unitlength}{4144sp}%
\begingroup\makeatletter\ifx\SetFigFont\undefined%
\gdef\SetFigFont#1#2#3#4#5{%
  \reset@font\fontsize{#1}{#2pt}%
  \fontfamily{#3}\fontseries{#4}\fontshape{#5}%
  \selectfont}%
\fi\endgroup%
\begin{picture}(1155,915)(256,-751)
\put(271,-691){\makebox(0,0)[lb]{\smash{{\SetFigFont{10}{12.0}{\familydefault}{\mddefault}{\updefault}{\color[rgb]{0,0,0}$I_1$}%
}}}}
\put(271, 29){\makebox(0,0)[lb]{\smash{{\SetFigFont{10}{12.0}{\familydefault}{\mddefault}{\updefault}{\color[rgb]{0,0,0}$I_4$}%
}}}}
\put(946, 29){\makebox(0,0)[lb]{\smash{{\SetFigFont{10}{12.0}{\familydefault}{\mddefault}{\updefault}{\color[rgb]{0,0,0}$I_3$}%
}}}}
\put(1396,-691){\makebox(0,0)[lb]{\smash{{\SetFigFont{10}{12.0}{\familydefault}{\mddefault}{\updefault}{\color[rgb]{0,0,0}$I_2$}%
}}}}
\put(766,-691){\makebox(0,0)[lb]{\smash{{\SetFigFont{10}{12.0}{\familydefault}{\mddefault}{\updefault}{\color[rgb]{0,0,0}$I_5,I_6$}%
}}}}
\end{picture}%
}
\end{center}

For any $\al=(1,1,1)$-generic $\te\in B$, the
fan $\Si_\te$ of $\lM_\te$ has five rays.
The matrix of $\ub\hi_\te^*:\cZ^{\Si_\te(1)}\arr M\dual$ equals
\smat{
0&2&1&0&1\\
0&0&1&1&0\\
1&1&1&1&1}
and is independent of the stability \te.
The Picard group $\Pic(\lM_\te)$ is isomorphic
to the cokernel of $\ub\hi_\te:M\arr\cZ^{\Si_\te(1)}$
(see \cite[Section 3.4]{Fulton2}).
We can choose a basis of $\Pic(\lM_\te)$ such
the matrix of $\cZ^{\Si_\te(1)}\arr\Pic(\lM_\te)$
is given by 
$$\smat{1 & 0 & 1 & -1 & -1\\ 0 & 1 & -1 & 1 & -1}.$$

There are $6$ different
chambers of \al-generic stabilities.
%(it is just a coincidence that the number of perfect matchings
%is the same).
Their representatives are given in Table \ref{spp:stab}.

\begin{table}[h!]
\begin{tabular}{|c|c|}
\hline
$N$ &$\te$\\
\hline
$1$&$(-2,1,1)$\\
$2$&$(1,-2,1)$\\
$3$&$(1,1,-2)$\\
$4$&$(2,-1,-1)$\\
$5$&$(-1,2,-1)$\\
$6$&$(-1,-1,2)$\\
\hline
\end{tabular}
\caption{\al-generic stabilities}
\label{spp:stab}
\end{table}

It is easy to see that the perfect matching
$I_5$ is stable with respect to 
$\te_2$, $\te_3$, and $\te_4$.
The perfect matching $I_6$ is stable with
respect to $\te_1$, $\te_5$, and $\te_6$.
The other perfect matchings are extremal
and therefore stable with respect to all $\te_i$,
$i=1,\dots,6$.

We will say that a pair of perfect matchings
is \te-stable if their union is \te-stable.
One can see that the pair $\set{I_1,I_3}$ is stable
only with respect to $\te_2$ and $\te_6$. The
pair $\set{I_2,I_4}$ is stable
only with respect to $\te_3$ and $\te_5$. This
uniquely determines the triangulation of the toric
diagram for any generic stability.
\begin{figure}[h!]%
\scalebox{1}{\input{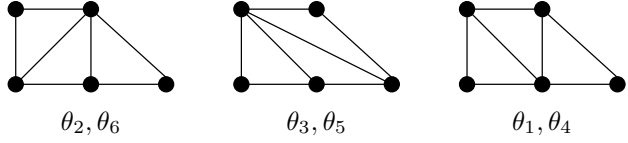}}
\caption{Triangulations of the toric diagram}%
\label{}%
\end{figure}

To construct the tilting bundle $\lU_\te$
on the moduli space $\lM_\te$, we choose paths
$u_1=e_1$ (trivial path), $u_2=12$, $u_3=13$ and intersect them
with \te-stable perfect matchings. This
gives us three vectors $\hi_\te(u_i)\in\cZ^{\Si_\te(1)}=\cZ^{5}$.
Their images with respect to the map
$\cZ^{\Si_\te(1)}\arr\Pic(\lM_\te)$ 
(given by the matrix $\smat{1 & 0 & 1 & -1 & -1\\ 0 & 1 & -1 & 1 & -1}$)
determine line bundles which are summands of $\lU_\te$.
The result is given by the following table

\begin{table}[h]
\begin{tabular}{|c|c|c|c|c|}
\hline
$\te$ &$\hi_\te(u_1)$&$\hi_\te(u_2)$&$\hi_\te(u_3)$&$\lU_\te$\\
\hline
$\te_2,\te_3,\te_4$
&$0$&$(1,0,0,0,1)$&$(0,1,0,0,1)$&$\lO\oplus\lO\smat{0\\-1}\oplus\lO\smat{-1\\0}$\\
$\te_1,\te_5,\te_6$
&$0$ &$(1,0,0,0,0)$&$(0,1,0,0,0)$&$\lO\oplus\lO\smat{1\\0}\oplus\lO\smat{0\\1}$\\
\hline
\end{tabular}
\end{table}

\subsection{Orbifold 
\texorpdfstring{$\cC^3/(\cZ_2\xx\cZ_2)$}{C3/Z2xZ2}}
Consider a finite abelian group $G=\cZ_2\xx\cZ_2$.
Consider an embedding $G\sb\SL_3(\cC)$, where the action
of the first copy of $\cZ_2$ is given by $\frac12(1,1,0)$
and the action of the second copy of $\cZ_2$ is given
by $\frac12(0,1,1)$. The character group $\hat G$ can be
identified with $\cZ_2\xx\cZ_2$.
We will use the following notation for the elements
of $\hat G$ (and sometimes for the elements of $G$):
$a=(1,0)$, $b=(0,1)$ and $c=(1,1)$.
The above representation of $G$ is isomorphic to $a\oplus b\oplus c$.

\begin{figure}[h!]
\scalebox{0.6}{\includegraphics{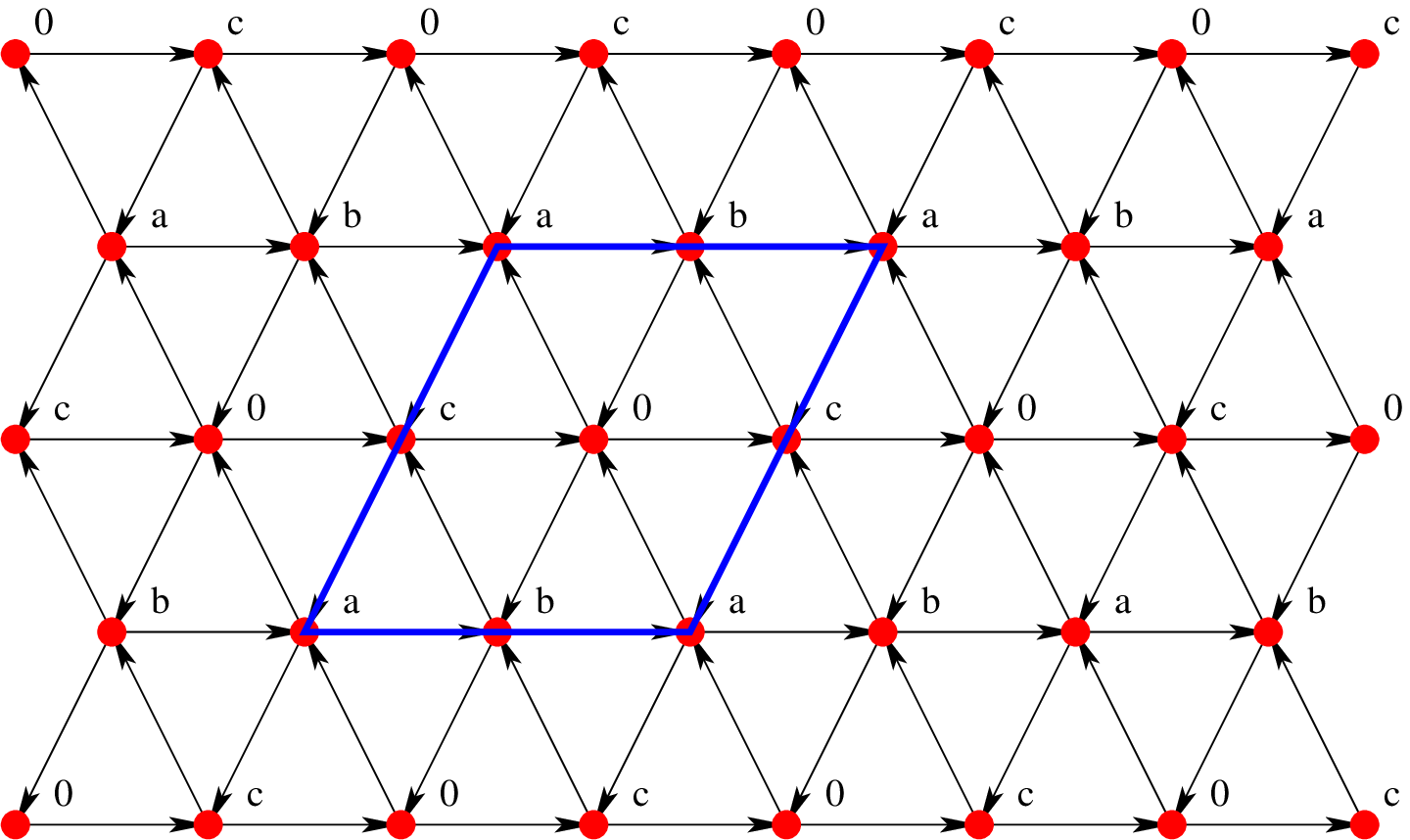}}
\caption{Periodic quiver for $\cC^3/(\cZ_2\xx\cZ_2)$}%
\label{fig:periodic1}%
\end{figure}

We can associate a brane tiling with the above embedding
(see \cite{Mozgov6} for the construction and notation).
The corresponding periodic quiver with a fundamental
domain is given in Figure \ref{fig:periodic1}.

The toric diagram of $\Spec\cC[M^+]$ is given in
Figure \ref{fig:toricd1} (see also \cite{Mozgov6}).
%We depicted there also the points
%of $M\dual\sb (M_0)_\cQ\dual$ contained in the triangle
%spanned by $e_1,e_2,e_3$.
\begin{figure}[h!]
\input{Z2xZ2.toric.pstex_t}
\caption{Toric diagram for $\cC^3/(\cZ_2\xx\cZ_2)$}%
\label{fig:toricd1}%
\end{figure}
The list of all possible perfect matchings
is given in Table \ref{table:pm2}.

\begin{table}[h!]
\begin{tabular}{|c|c|c|c|}
\hline
$N$ &$I$& $4\ub\hi_I$ & \\
\hline
$1$ & $cb,0a,0b,ca$ & $(2,2,0)$ &$a$\\%a
$2$ & $cb,0a,bc,a0$ & $(4,0,0)$ &$e_1$\\
$3$ & $cb,ab,c0,a0$ & $(2,0,2)$ &$c$\\%c+
$4$ & $ac,b0,0b,ca$ & $(0,4,0)$ &$e_2$\\
$5$ & $ac,b0,bc,a0$ & $(2,2,0)$ &$a$\\%a+
$6$ & $ac,ab,0b,0c$ & $(0,2,2)$ &$b$\\%b
$7$ & $ba,b0,c0,ca$ & $(0,2,2)$ &$b$\\%b+
$8$ & $ba,0a,bc,0c$ & $(2,0,2)$ &$c$\\%c
$9$ & $ba,ab,c0,0c$ & $(0,0,4)$ &$e_3$\\
\hline
\end{tabular}
\caption{Perfect matchings for $\cC^3/(\cZ_2\xx\cZ_2)$}
\label{table:pm2}
\end{table}

We will consider only stabilities 
$$\te_1=(-3,1,1,1),\quad \te_2=(-3,-1,2,2),\quad 
\te_3=(-2,3,1,-2),$$
where the order of the coordinates of $\cZ^{Q_0}=\cZ^{\hat G}$ is
given by $0,a,b,c$.
Note that $\lM_{\te_1}$ is isomorphic to
$\Hilb^G(\cC^3)$.

A subset $I\sb Q_1$ is $\te_1$-stable
if and only if there exists a path in $Q\ms I$ from vertex $0\in Q_0$ to
any other vertex of $Q$.
The non-extremal $\te_1$-stable perfect matchings
are $I_3,I_5,I_7$.
A subset $I\sb Q_1$ is $\te_2$-stable
if and only if there exist paths in $Q\ms I$
from $0$ to $b$ and $c$ and from $a$
to $b$ or $c$.
The non-extremal $\te_2$-stable perfect matchings are
the same as for $\te_1$.
A subset $I\sb Q_1$ is $\te_3$-stable
if and only if there exist paths in $Q\ms I$
from $0$ to $a$, from $c$ to $a$, and
from $0$ or $c$ to $b$.
The non-extremal $\te_3$-stable perfect matchings
are $I_3,I_5,I_6$.

\begin{figure}[h!]
\input{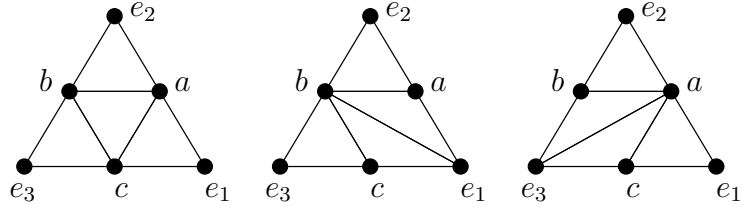}
\caption{Toric diagrams for $\te_1,\te_2,\te_3$}%
\label{fig:toricdx}%
\end{figure}

To determine the toric diagram of $\lM_{\te_i}$, $i=1,2,3$,
we have to find such pairs of $\te_i$-stable perfect
matchings that their union is still $\te_i$-stable (we call such
pairs $\te_i$-stable).
For $\te_1$, the pairs
$\set{I_3,I_5},\set{I_5,I_7},\set{I_3,I_7}$
are stable.
The toric diagram of $\lM_{\te_1}$
is given in Figure \ref{fig:toricdx}.
For $\te_2$, the pair 
$\set{I_3,I_5}$ (corresponds to the edge $ca$) is non-stable.
This uniquely determines the toric diagram of $\lM_{\te_2}$
(see Figure \ref{fig:toricdx}).
%Note that the toric diagrams of $\lM_{\te_1}$ and
%$\lM_{\te_2}$ are different, also the sets
%of $\te_1$-stable perfect matchings and
%$\te_2$-stable perfect matchings coincide. 
For $\te_3$, the pair $\set{I_3,I_6}$ (corresponds to the
edge $cb$) is non-stable.
This uniquely determines the toric diagram of $\lM_{\te_3}$
(see Figure \ref{fig:toricdx}).

\begin{figure}[h]
\scalebox{0.8}{\includegraphics{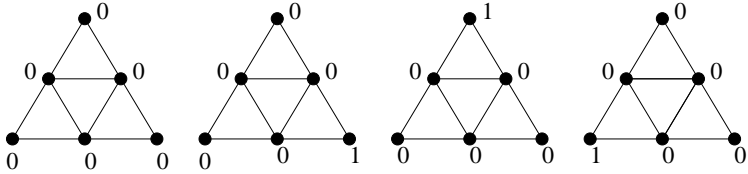}}
\caption{Cartier divisors for a tilting collection on $\lM_{\te_1}$}%
\label{fig:cartier1}%
\end{figure}

\begin{figure}[h]
\scalebox{0.8}{\includegraphics{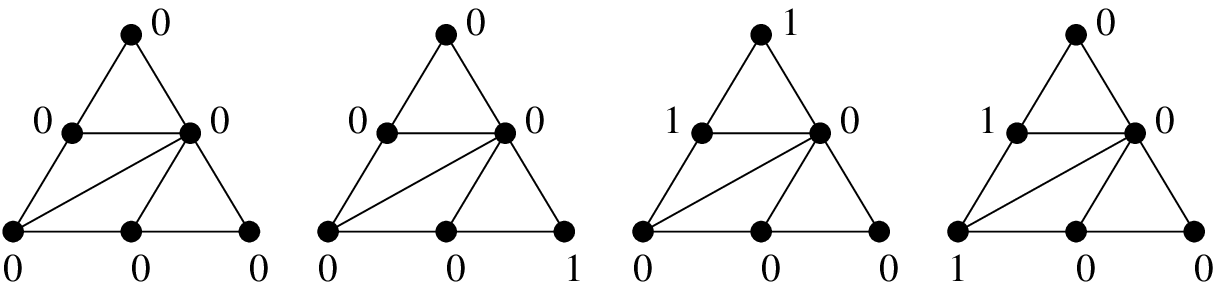}}
\caption{Cartier divisors for a tilting collection on $\lM_{\te_3}$}%
\label{fig:cartier2}%
\end{figure}

To determine the tilting bundle $\lU_{\te_i}$ on $\lM_{\te_i}$,
$i=1,2,3$, we choose
paths from vertex $0\in Q_0$ to all other vertices
of $Q$ and intersect these paths with $\te_i$-stable
perfect matchings. We choose paths
$e_0,0a,0b,0c$. The result of intersecting
these paths with $\te_1$-stable perfect matchings
is given in Figure \ref{fig:cartier1}.
In this way we get Cartier divisors for a tilting collection on $\lM_{\te_1}$.
The result for $\te_2$ is the same, as
$\te_1$-stable perfect matchings and
$\te_2$-stable perfect matchings coincide.
The result for $\te_3$ is given in Figure \ref{fig:cartier2}.
This gives us Cartier divisors for a tilting collection on $\lM_{\te_3}$.

 %examples

\bibliography{../tex/fullbib}
\bibliographystyle{../tex/hamsplain}
\end{document}